\documentclass[12pt]{amsart}
\usepackage{amssymb}
\usepackage{float}
\usepackage{amsmath}
\usepackage{graphicx}
\usepackage{amscd}

\theoremstyle{plain}
\newtheorem{X}{X}[section]

\newtheorem{conjecture}[X]{Conjecture}

\theoremstyle{definition}

\newtheorem{remark}[X]{Remark}

\numberwithin{equation}{section}
\numberwithin{equation}{section}
\renewcommand{\tfrac}{\frac}
\usepackage{verbatim}

\makeindex

\newenvironment{mylist}
{\begin{list}
{--}
{\setlength{\leftmargin}{.5in}\setlength{\rightmargin}{.5in}}}
{\end{list}}

\renewenvironment{itemize}
{\begin{mylist}}
{\end{mylist}}

\numberwithin{equation}{section}

\renewcommand{\frac}{\tfrac}

\def\1{{1\mkern-7mu1}}  

\newcommand\Aut{\operatorname{Aut}}

\newcommand\End{\operatorname{End}}

\newcommand\Gal{\operatorname{Gal}}
\newcommand\GL{\operatorname{GL}}

\newcommand\Hom{\operatorname{Hom}}

\newcommand\Rep{\operatorname{Rep}}

\newcommand\Sh{\operatorname{Sh}}

\def\CM{{\bold{CM}}}

\def\LCM{{\bold{LCM}}}

\def\LMot{{\bold{LMot}}}

\def\Mot{{\bold{Mot}}}

\def\Rep{{\bold{Rep}}}
\def\Vc{{\bold{Vec}}}

\begin{document}
\title[ihp00 April 30, 2000]{Towards a proof of the conjecture of Langlands and Rapoport}
\dedicatory{April 30, 2000; Rough first draft.}
\author[ihp00 April 30, 2000]{J.S. Milne}
\address{2679 Bedford Rd., Ann Arbor, MI 48104, USA.}
\email{math@jmilne.org, \emph{Web site:}\texttt{www.jmilne.org/math/}}
\thanks{Thanks ... \\
Text for a talk April 28, 2000, at the Conference on Galois Representations,
Automorphic Representations and Shimura Varieties, Institut Henri
Poincar\'e, Paris, April 24--29, 2000.}
\maketitle
\tableofcontents

\section{Introduction}

A reductive group $G$ over $\mathbb{Q}{}$ plus a $G(\mathbb{R}{})$-conjugacy
class $X$ of homomorphisms $\mathbb{C}{}^{\times }\rightarrow G(\mathbb{R}%
{}) $ satisfying certain axioms (Deligne 1979) defines a Shimura variety $\Sh%
(G,X)$, which is the projective system of the double coset spaces 
\begin{equation*}
\quad \Sh_{K}(G,X)=G(\mathbb{Q}{})\backslash X\times G(\mathbb{A}_{f})/K,
\end{equation*}

\noindent with $K$ running over the compact open subgroups of $G(\mathbb{A}%
_{f})$. \noindent The axioms imply that the $X$ has the structure of a
disjoint union of bounded symmetric domains, and \noindent $G(\mathbb{Q}%
{})\backslash X\times G(\mathbb{A}_{f})/K$ is a disjoint union of spaces of
the form $\Gamma \backslash X^{+}$ with $X^{+}$ a connected component of $X$
and $\Gamma $ a congruence subgroup of $G^{\text{der}}(\mathbb{R)}$, and so $%
\Sh(G,X)$ is a projective system of analytic spaces. The theorem of Baily
and Borel shows that $\Sh(G,X)$ is a projective system of algebraic
varieties (not connected!) over $\mathbb{C}{}$. A theorem of Shimura,
Deligne, et al. shows that $\Sh(G,X)$ has a canonical model over a certain
number field $E(G,X)$, called the reflex field. Thus, we may reduce $\Sh%
(G,X) $ modulo a prime ideal $\mathfrak{p}_{v}$ of $E(G,X)$, for
example, by embedding each $\Sh_{K}(G,X)$ in projective space and then
scaling and reducing the equations modulo $\mathfrak{p}_{v}$, to obtain
a projective family of varieties over the residue field $\kappa (v)$.

However, without any further conditions on $G$ and $X$, the reduced
varieties may be very singular. To avoid this, we assume that a hyperspecial
group $K_{p}$ has been given, and consider the system 
\begin{equation*}
\Sh_{p}(G,X)=\{\Sh_{K^{p}\cdot K_{p}}(G,X)\mid K^{p}\text{ compact open in }%
G(\mathbb{A}_{f}^{p})\}.
\end{equation*}

\noindent Here $\mathbb{A}_{f}^{p}$ is the ring of finite ad\`{e}les with
the $p$-component omitted. The existence of $K_{p}$ implies that $G$ is
unramified over $\mathbb{Q}{}_{p}$, and by considering only groups of the
form $K^{p}\cdot K_{p}$ we are, in effect, imposing a level structure only
away from $p$.

When we reduce modulo a prime $\mathfrak{p}_{v}$ dividing $p$, we obtain
pro-variety $\Sh_{p}(G,X)_{v}$ over $\kappa (v)$. We write $\Sh_{p}(\mathbb{F%
}{})$ for the set of its points in the algebraic closure $\mathbb{F}{}$ of $%
\kappa (v)$. This is a set with an action of the Frobenius generator of $\Gal%
(\mathbb{F}{}/\kappa (v))$ and $G(\mathbb{A}_{f}^{p})$. The conjecture of
Langlands and Rapoport gives a description $\Sh_{p}(\mathbb{F}{})$, together
with the two actions, directly in terms of the initial data $G,X,K_{p}.$

\section{The Different Types of Shimura Varieties}

We shall look at the conjecture in three cases. Generally, I shall regard
abelian varieties as lying in the category of abelian varieties up to
isogeny, i.e., the category whose objects are the abelian varieties but in
which the Hom-sets have been tensored with $\mathbb{Q}{}$.

\subsubsection{Shimura varieties of PEL-type.}

For an appropriate choice of a representation $G\hookrightarrow \GL(V)$,
Shimura varieties of PEL-type become moduli varieties\footnote{%
Strictly, I should say pro-variety...} in characteristic zero, namely, their
points classify isomorphism classes of triples $(A,\eta ^{p},\Lambda _{p})$
where $A$ is an abelian variety endowed with a polarization and an action of
a fixed $\mathbb{Q}{}$-algebra $B$, $\eta ^{p}$ is a prime-to-$p$ level
structure on $A$, and $\Lambda _{p}$ is a lattice in $H_{1}(A_{\text{et}},%
\mathbb{Q}{}_{p})$. The triples are required to satisfy certain conditions,
for example, that the representation of the $\mathbb{Q}{}$-algebra $B$ on
the tangent space to $A$ at zero lies in a fixed isomorphism class. An
isomorphism of triples is an isogeny of abelian varieties (element of $\Hom%
(A,A^{\prime })\otimes \mathbb{Q}{}$ with an inverse in $\Hom(A^{\prime
},A)\otimes \mathbb{Q}{}$) preserving all the structure.

\subsubsection{Shimura varieties of Hodge type.}

This class has description similar to the preceding class, except that now
(in characteristic zero) $\Sh_{p}(G,X)$ is the moduli variety for triples $%
(A,\eta ^{p},\Lambda _{p})$ where $A$ is an abelian variety endowed with
some Hodge classes (in the sense of Deligne 1982).

\subsubsection{Shimura varieties of abelian type.}

This is the almost-general case, since it excludes only the Shimura
varieties defined by groups with factors of type $E_{6}$, $E_{7}$, and
certain mixed types $D$. Associated with the datum defining the Shimura
variety there is a ``weight'' homomorphism $w_{X}\colon \mathbb{G}%
_{m}\rightarrow G$. When $w_{X}$ is defined over $\mathbb{Q}{}$, then (in
characteristic zero) the choice of a representation $G\hookrightarrow \GL(V)$
realizes $\Sh_{p}(G,X)$ as the moduli variety for triples $(M,\eta
^{p},\Lambda _{p})$ where $M$ is now a motive rather than an abelian variety
(Milne 1994). When $w_{X}$ is not defined over $\mathbb{Q}{}$, then $%
Sh_{p}(G,X)$ is not a moduli variety.

\subsubsection{Comments.}

There are many more Shimura varieties of abelian type than of Hodge type,
and many more of Hodge type than of PEL-type. The case of Hodge type is
always a useful stepping stone to the almost-general case. The PEL case is
included only because it is so much easier than the other cases, and
studying it gives a guide to how to proceed in the other cases.

Surprisingly, the study of Shimura varieties of Hodge type turns out to be
much more difficult than that of those PEL-type for two reasons. First,
multilinear algebra is more difficult than linear algebra. For example, as I
mentioned above, one of the conditions for a triple $(A,\eta ^{p},\Lambda
_{p})$ to lie in the family parametrized by $\Sh_{p}(G,X)$ is that the
representation of the algebra $B$ on the tangent space be of a fixed
isomorphism type. Stating such a condition for tensors in some spaces is
more difficult. The second reason is that Deligne's Hodge classes are
defined only in characteristic zero.

\section{The statement of the conjecture of Langlands and Rapoport}

Langlands and Rapoport first ``define'' a category of motives $\Mot(\mathbb{F%
})$ over $\mathbb{F}{}$ --- the reason for the quotes will be explained
below.

Assume that the weight $w_{X}$ for $\Sh_{p}(G,X)$ is defined over $\mathbb{Q}%
{}$. Then, the above discussion suggests that, if we fix a representation $%
G\hookrightarrow \GL(V)$, then there should be a one-to-one correspondence
the points in $\Sh_{p}(\mathbb{F}{})$ and a set of isomorphism classes of
triples $(M,\eta ^{p},\Lambda _{p})$ with $M$ a motive over $\mathbb{F}{}$, $%
\eta ^{p}$ a prime-to-$p$ level structure on $M$, and $\Lambda _{p}$ a
lattice in $H_{\text{crys}}^{\ast }(M)$.

When, we vary the representation of $G$, then this should become a
one-to-one correspondence between $\Sh_{p}(\mathbb{F)}$ and a set of
isomorphism classes of triples $($\underline{$M$}$,\eta ^{p},\Lambda _{p})$
where \underline{$M$} is now a functor $\Rep(G)\rightarrow \Mot(\mathbb{F})$.

When we choose a fibre functor $\omega $ for $\Mot(\mathbb{F})$ and let $%
\mathfrak{P}$ be the corresponding groupoid, $\mathfrak{P}{}=_{\text{df}}%
\Aut_{\mathbb{Q}{}}^{\otimes }(\omega )$, then the theory of Tannakian
categories shows that to give an \underline{$M$} is the same as to give a
morphism of groupoids $\phi \colon \mathfrak{P}{}\rightarrow \mathfrak{G}{}_{G}$.

Finally, Langlands and Rapoport define another groupoid $\mathfrak{Q}{}$ having $%
\mathfrak{P}{}$ as a quotient, to allow for the weight to be irrational. They
define a set of triples $(\phi ,\eta ^{p},\Lambda _{p})$ (depending only on $%
G,X,K_{p}$) where $\phi $ is now a homomorphism $\mathfrak{Q}{}\rightarrow \mathfrak{%
\mathfrak{G}{}}_{G}$ and their conjecture states that the elements of $\Sh_{p}(%
\mathbb{F}{})$ should be in one-to-one correspondence with the set $LR(%
\mathbb{F}{})$ of isomorphism classes of these triples. There is a natural
action of the Frobenius automorphism and of $G(\mathbb{A}_{f}^{p})$ on the
triples, and the correspondence should respect these actions.

\begin{remark}
Of course, it is easy to guess that somehow a Shimura variety modulo a prime 
$\mathfrak{p}{}_{v}$ should parametrize isomorphism classes of motives with
additional structure. The point of the paper of Langlands and Rapoport is to
define a category of motives $\Mot(\mathbb{F})$ and then to state precise
conditions on the triples $(\phi ,\eta ^{p},\Lambda _{p})$ that are to
correspond to a given Shimura variety.
\end{remark}

\begin{remark}
I will discuss the definition of $\Mot(\mathbb{F})$ below. The statement of
the precise conditions on the triples $(\phi ,\eta ^{p},\Lambda _{p})$ is
quite complicated, especially that on $\Lambda _{p}$, and I refer the reader
to the original paper or Milne 1992 for these. Here I will only make a few
comments.

Langlands and Rapoport defined the notion of an \emph{admissible }
homomorphism $\phi \colon \mathfrak{Q}{}\rightarrow \mathfrak{G}{}_{G}$. To be
admissible, a homomorphism must satisfy one condition for each prime $l$
(including $l=\infty $) and a condition on the composite of $\phi $ with $%
\mathfrak{G}{}_{G}\rightarrow \mathfrak{G}_{G/G^{\text{der}}}$.\newline

A special point in $X$ defines a homomorphism $\phi \colon \mathfrak{Q}%
{}\rightarrow \mathfrak{G}_{G}$, called \emph{special}. Langlands and
Rapoport show that, when $G^{\text{der}}$ is simply connected, a
homomorphism is admissible if and only if it is isomorphic to a special
homomorphism. They also show that if their conjecture is correct for groups
with $G^{\text{der}}$ simply connected, then it can't be correct without
this condition (so they stated their conjecture only for $G^{\text{der}}$
simply connected).

I showed, on the other hand, that if one replaces ``admissible'' with the
condition ``is isomorphic to a special homomorphism'' then the conjecture is
true for all Shimura varieties when it is true for those with $G^{\text{der}}
$ simply connected. This allows us to state the conjecture for all Shimura
varieties and (to some extent) reduces the problem of its proof to the case
where $G^{\text{der}}$ is simply connected.
\end{remark}

\begin{remark}
The description given by the conjecture of Langlands and Rapoport may look
too abstract and complicated to be of use, but, in fact, it fairly
straightforward to derive formulas for the numbers of points in terms of
orbital and twisted orbital integrals (as conjectured by Langlands and
Kottwitz) from it. This is explained in Milne 1992.

\emph{From now on, I will assume that }$G^{\text{der}}$ \emph{simply
connected.}
\end{remark}

\section{Improvements to the Statement of the Conjecture}

\subsubsection{Canonical integral models}

One defect of the original conjecture is that it doesn't specify how to
reduce the Shimura variety. Defining a reduction amounts to defining a model
of $\Sh_{p}(G,X)$ over the ring of integers in $E_{v}$. Evidently, if the
conjecture is true for one integral model, it will be false for most others.
I suggested (Milne 1992) that there should be a canonical integral model
characterized by a certain N\'{e}ron-type property. The existence of such a
model has been proved by Vasiu for $p\geq 5$ (Vasiu 1999).

To the original conjecture, one should add that the reduction is that
defined by the canonical integral model.

\subsubsection{The definition of $\Mot(\mathbb{F})$}

The Tate conjecture implies that the category of motives over $\mathbb{F}{}$
should be Tannakian with a certain specific protorus $P$ as its band. The
Tannakian categories with $P$ as band are classified up to $P$-equivalence
by the cohomology group $H^{2}(\mathbb{Q}{},P)$, and Langlands and Rapoport
showed that there is only one class in $H^{2}(\mathbb{Q}{},P)$ giving a
Tannakian category for which the correct fibre functors exist. They define $%
\Mot(\mathbb{F})$ to be any Tannakian category with band $P$ having this
cohomology class. Thus, $\Mot(\mathbb{F})$ is only defined up to a nonunique 
$P$-equivalence. Because some $H^{1}$'s vanish, the category is a little
better defined than one might expect but it still not possible to talk of
objects in $\Mot(\mathbb{F})$.

For example, let $M_{1}$ be one model for $\Mot(\mathbb{F})$ and let $X$ be
an object of $M_{1}$. If $M_{2}$ is a second model, then there is a $P$%
-equivalence $F\colon M_{1}\rightarrow M_{2}$, and so $X$ corresponds to an
object $FX$ in $M_{2}$. But, there is no special $F$, and if $F^{\prime }$
is second $P$-equivalence $M_{1}\rightarrow M_{2}$, then $X$ will correspond
to a second object $F^{\prime }X$ of $M_{2}$. The objects $FX$ and $%
F^{\prime }X$ will be isomorphic, but there is no preferred isomorphism.
Thus, all one say is that isomorphism classes of objects in $M_{1}$
correspond to isomorphism classes in $M_{2}$. This is scarcely better than
the information provided by the Honda-Tate theorem on the category of
abelian varieties up to isogeny over $\mathbb{F}{}$: it classifies only the
isomorphism classes and their endomorphism algebras.

Below I shall provide a more precise definition of $\Mot(\mathbb{F})$ for
which, in the above discussion, $FX$ and $F^{\prime }X$ will isomorphic with
a \emph{unique }isomorphism. In other words, the object in $M_{2}$
corresponding to $X$ in $M_{1}$ will be well-defined up to a unique
isomorphism.

This more precisely defined category (which gives a $\mathfrak{P}{}$ and $\mathfrak{Q%
}{}$) should be the one used in the statement of the conjecture.

\subsubsection{Canonicalness}

Once one has the notion of a canonical integral model, the set with
operators $\Sh_{p}(\mathbb{F}{})$ is \emph{canonically }associated with $%
G,X,K_{p}$. Moreover, once one chooses a fibre functor on $\Mot(\mathbb{F})$%
, the set $LR(\mathbb{F}{})$ is also \emph{canonically }associated with $%
G,X,K_{p}$. Clearly, one should require that the one-to-one correspondence
in the statement of the conjecture of Langlands and Rapoport \emph{be
canonical}.

In fact, I expect one can prove a uniqueness statement of the following
form: there is at most one family of bijections $LR(G,X,K_{p})(\mathbb{F}%
{})\rightarrow Sh_{p}(G,X)(\mathbb{F}{})$ having certain functoriality
properties and giving the correct map for the Siegel modular variety.

\bigskip

Henceforth, by the conjecture of Langlands and Rapoport I mean the canonical
conjecture.

\section{The Conditional Proof of Langlands and Rapoport in the PEL-case}

In their paper (1987), Langlands and Rapoport prove their conjecture for
Shimura varieties of PEL-type under the assumption of:

\begin{enumerate}
\item  the Hodge conjecture for complex abelian varieties of CM-type;

\item  the Tate conjecture for abelian varieties over $\mathbb{F}{}$;

\item  Grothendieck's standard conjectures for abelian varieties over $%
\mathbb{F}{}$.
\end{enumerate}

One of the difficulties is that their abstractly-defined category $\Mot(%
\mathbb{F})$ does not contain the category of abelian varieties (up to
isogeny) in any natural way. When one assumes (b), one gets a well-defined
category of motives containing the category abelian varieties (up to
isogeny), namely the category of motives based on abelian varieties using
the algebraic classes modulo numerical equivalence as correspondences.

Another difficulty is that Deligne's Hodge classes make sense only in
characteristic zero. Let $\CM(\mathbb{Q}{}^{\text{al}})$ be the category of
motives based on abelian varieties over $\mathbb{Q}^{\text{al}}$ of CM-type
using the Hodge classes as correspondences. When (a) is assumed, the Hodge
classes will be algebraic, and therefore will reduce to algebraic classes.
We then get a canonical functor $R\colon \CM(\mathbb{Q}{}^{\text{al}%
})\rightarrow \Mot(\mathbb{F})$.

In summary, assuming (a) and (b) we get a canonical commutative diagram 
\begin{equation*}
\begin{array}{ccc}
\CM(\mathbb{Q}{}^{\text{al}}) & \overset{I}{\leftarrow } & \LCM(\mathbb{Q}%
{}^{\text{al}}) \\ 
\downarrow R &  & \downarrow R \\ 
\Mot(\mathbb{F}) & \overset{I}{\leftarrow } & \LMot(\mathbb{F})
\end{array}
\end{equation*}

\noindent where $\LCM(\mathbb{Q}{}^{\text{al}})$ is the category of motives
based on abelian varieties of CM-type over $\mathbb{Q}{}^{\text{al}}$ using
the Lefschetz classes\footnote{%
A Lefschetz class is an element of the $\mathbb{Q}{}$-algebra generated by
divisor classes inside the $\mathbb{Q}{}$-algebra of algebraic classes
modulo numerical equivalence (or inside a Weil cohomology --- there is no
difference for abelian varieties).} as correspondences, and $\LMot(\mathbb{F}%
)$ is the similar category based on abelian varieties over $\mathbb{F}{}$.

Finally, recall that a Weil form on an object $X$ of a Tannakian category is
a form (bilinear or sesquilinear according to context) that induces a
positive involution on $\End(X)$, and that to give a polarization on a
Tannakian category is to give a distinguished class of ``positive'' Weil
forms for each object satisfying certain compatibility conditions. A
polarization of an abelian variety $A$ in the sense of algebraic geometry
defines a Weil form on $h_{1}A$. Grothendieck's standard conjectures imply
that there is a polarization on $\Mot(\mathbb{F})$ for which these geometric
Weil forms are all positive.

With the assumption of (a), (b), (c), the proof of the conjecture of
Langlands and Rapoport for Shimura varieties of PEL-type becomes fairly
straightforward. The canonical integral model is, in this case, a moduli
variety, and so a point in $\Sh_{p}(\mathbb{F}{})$ corresponds to an
isomorphism class of triples $(A,\lambda ^{p},\Lambda _{p})$ with $A$ an
abelian variety endowed with a polarization and an action of a $\mathbb{Q}{}$%
-algebra $B$. The object $I(h_{1}A)$ then defines (by the theory of
Tannakian categories) a morphism $\phi _{A}\colon \mathfrak{P}{}\rightarrow 
\mathfrak{G}{}_{G}$, and Langlands and Rapoport verify that there is a canonical
bijection between the set $S(A)$ of isomorphism classes of triples $%
(A,\lambda ^{p},\Lambda _{p})$ (fixed $A$) and the set of isomorphism
classes of triples $(\phi _{A},\lambda ^{p},\Lambda _{p})$ (fixed $\phi _{A}$%
). The abelian variety $A$ with its PE-structure is, almost by definition,
the reduction (up to isogeny) of an abelian variety $\tilde{A}$ with
PE-structure in the family parametrized by $\Sh_{p}$ in characteristic zero.
A theorem of Zink's shows that $\tilde{A}$ can be chosen to be of CM-type.
It therefore corresponds to a special point $x$ of $X$, and one verifies
that $\phi _{x}\approx \phi _{A}$. Thus $\phi _{A}$ is admissible.

\section{Towards an unconditional proof in the PEL-case}

We wish to carry out the above argument without assuming (a), (b), or (c).
In Milne 1999, I showed that (a) implies (b), and I can show that (a)
implies at least the consequence of (c) needed for the above proof. Thus,
instead of three conjectures we need to assume only one, namely, the Hodge
conjecture for abelian varieties of CM-type. Unfortunately, the meagre
progress made on the Hodge conjecture in the 50 years since the conjecture
was made suggests that it will not be wise to wait for a proof of the Hodge
conjecture, even for abelian varieties of CM-type.

First I explain my new construction of $\Mot(\mathbb{F})$. The conjecture
(a) implies that we have a functor $R\colon \CM(\mathbb{Q}{}^{\text{al}%
})\rightarrow \Mot(\mathbb{F})$ bound by a map $P\hookrightarrow S$ of
pro-tori. The group $P$ acts on the objects of $\Mot(\mathbb{F})$ and we let 
$\Mot(\mathbb{F})^{P}$ be the subcategory of objects on which $P$ acts
trivially. Thus $\Mot(\mathbb{F})^{P}$ comprises the motives consisting
entirely of algebraic classes. Let $\1$ be an identity object of $\Mot(%
\mathbb{F})$. Then $X\mapsto \Hom(\1,X)$ is a fibre functor $\Mot(\mathbb{F}%
)^{P}\rightarrow \Vc_{\mathbb{Q}{}}$. Its composite with $R$ is a fibre
functor $\omega _{0}$ on $\CM(\mathbb{Q}{}^{\text{al}})$, and I claim we can
(essentially) reverse this procedure and reconstruct $\Mot(\mathbb{F})$ from 
$\CM(\mathbb{Q}{}^{\text{al}})$ and $\omega _{0}$. I now drop all
assumptions.

First, one shows that there is a fibre functor $\omega _{0}$, unique up to
isomorphism, that when tensored with $\mathbb{Q}{}_{l}$ is in the
``correct'' isomorphism class for all $l\leq \infty $. Now define $\Mot(%
\mathbb{F})^{\prime }$ as follows:

\begin{itemize}
\item  $\Mot(\mathbb{F})^{\prime }$ has one object $\bar{X}$ for each object 
$X$ of $\Mot(\mathbb{F})$;

\item  for objects $\bar{X}$, $\bar{Y}$ of $\Mot(\mathbb{F})$, define $\Hom(%
\bar{X},\bar{Y})=\omega _{0}($\underline{$\Hom$}$(X,Y)^{P})$.
\end{itemize}

\noindent Here \underline{$\Hom$}$(X,Y)$ is the internal Hom of $X$ and $Y$
in $\CM(\mathbb{Q}{}^{\text{al}})$, and \underline{$\Hom$}$(X,Y)^{P}$ is the
largest subobject fixed by $P$. Now $\Mot(\mathbb{F})$ is obtained from $\Mot%
(\mathbb{F})^{\prime }$ by adding the images of projectors, i.e., by taking
the pseudo-abelian (Karoubian) hull. It is only an exercise, using the
dictionary between Tannakian categories and gerbs, to show that this does
gives a Tannakian category and that $X\mapsto \bar{X}$ defines a tensor
functor $R\colon \CM(\mathbb{Q}{}^{\text{al}})\rightarrow \Mot(\mathbb{F})$
bound by $P\hookrightarrow S$.

Because $\omega _{0}$ is uniquely determined up to an isomorphism, which
itself is determined up to a \emph{unique }isomorphism, $\Mot(\mathbb{F})$
has the uniqueness property claimed.

To give a fibre functor on $\Mot(\mathbb{F})$ is to give a fibre functor $%
\omega $ on $\CM(\mathbb{Q}{}^{\text{al}})$ together with an isomorphism $%
\omega _{0}\rightarrow \omega \circ R$. In this way, one obtains fibre
functors $\omega _{\ell }\colon \Mot(\mathbb{F})\rightarrow \Vc_{\mathbb{Q}%
{}_{\ell }}$ for each $\ell \neq p$, and a functor $\omega _{p}\colon \Mot(%
\mathbb{F})\rightarrow Isoc(\mathbb{F}{})$, well-defined to isomorphism.

We thus have: 
\begin{equation*}
\begin{array}{ccc}
\CM(\mathbb{Q}{}^{\text{al}}) & \overset{I}{\leftarrow } & \LCM(\mathbb{Q}%
{}^{\text{al}}) \\ 
\downarrow R &  & \downarrow R \\ 
\Mot(\mathbb{F}) &  & \LMot(\mathbb{F})
\end{array}
\end{equation*}

\noindent The subcategory of $\Mot(\mathbb{F})$ of objects of weight $1$ is
certainly equivalent to the category of abelian varieties up to isogeny over 
$\mathbb{F}$, and it follows one does get a functor $I\colon \LMot(\mathbb{F}%
)\rightarrow \Mot(\mathbb{F})$. One can even choose it so that $\omega
_{l}^{M}\circ I\approx \omega _{l}^{LM}$ for all $l$. However, without
something extra, $I$ will not be canonical and the diagram may not quite
commute (its failure to commute is measured by a class in $H^{1}(\mathbb{Q}%
{},T)$, $T$ the fundamental group of $\LCM(\mathbb{Q}{}^{\text{al}})$, that
is trivial at all the finite primes).

On applying the method of Langlands and Rapoport described above, one
obtains a description of $\Sh_{p}(\mathbb{F}{})$ as the set of isomorphism
classes of triples $(\phi ,\eta ^{p},\Lambda _{p})$ exactly as conjectured,
except that the $\phi $ need not be admissible (each $\phi $ may be a twist
of an admissible $\phi $ by a cohomology class which may be chosen to come
from the centre of $I_{\phi }$ and split at all the finite primes). Thus,
one doesn't obtain the canonical LR conjecture. For that, one needs the
following conjecture.

\begin{conjecture}[A]
Let $A$ be an abelian variety of CM-type over $\mathbb{Q}{}^{\text{al}}$
(say), and let $\alpha $ be a Hodge class on $A$ (thus $\alpha $ lies in a
certain $\mathbb{Q}{}$-vector space). Such an $\alpha $ defines cohomology
classes $\alpha _{l}$ for all $l$ (in $H^{2\ast }(A_{\text{et}},\mathbb{Q}%
_{\ell }(\ast ))$ for $\ell \neq p$ and in $H_{\text{dR}}^{2\ast }(A)$ for $%
\ell =p$). Each $\alpha _{l}$ defines a cohomology class $\bar{\alpha}_{l}$
on the reduction $A_{\mathbb{F}{}}$ of $A$ ($\bar{\alpha}_{p}$ lies in the
crystalline cohomology). Let $\alpha _{f}=(\alpha _{l})$ and $\bar{\alpha}%
_{f}=(\bar{\alpha}_{l}).$

If $\bar{\alpha}_{f}$ is in $\mathbb{A}{}_{f}$-span of the Lefschetz classes
on $A_{\mathbb{F}{}}$, then it is a Lefschetz class (i.e., it is in the $%
\mathbb{Q}{}$-space of such classes).

Equivalently, if $\bar{\alpha}_{l}$ is in the $\mathbb{Q}{}_{l}$-span of the
Lefschetz classes for all $l$, then $\bar{\alpha}_{l}$ is the cohomology
class of a Lefschetz class for all $l$, which is independent of $l$.
\end{conjecture}

Roughly speaking, the conjecture says that Hodge classes on $A$ that look as
though they should become Lefschetz on $A_{\mathbb{F}{}}$, do in fact become
Lefschetz. It is a compatibility conjecture between a $\mathbb{Q}{}$%
-structure in characteristic zero and a $\mathbb{Q}{}$-structure in
characteristic $p$.

When Conjecture A is assumed, $\Mot(\mathbb{F})$ becomes well-defined up to
a unique equivalence, the diagram commutes, there are well-defined functors $%
\omega _{l}$ on $\Mot(\mathbb{F})$ composing correctly with $R$ and $I$, and
there is a unique polarization on $\Mot(\mathbb{F})$ for which the Weil
forms coming from algebraic geometry are positive. Thus, the situation is
essentially as good as when one assumes (a), (b), (c), and the argument in
Langlands and Rapoport does give a proof of the canonical form of their
conjecture for Shimura varieties of PEL-type.

I hope that Conjecture A is susceptible to proof by the same methods that
Deligne used to prove his result on Hodge classes. Specifically, I can show
that it is true when $\alpha $ is algebraic (and so the conjecture is
implied by the Hodge conjecture for abelian varieties of CM-type). Moreover,
I can show that there is a subgroup $G$ of the Lefschetz group of $A$ such
that a Hodge class $\alpha $ becomes ($\mathbb{Q}{}$-rationally) Lefschetz
on $A_{\mathbb{F}{}}$ if and only if $\alpha $ is fixed by $G$. The next
step will be to prove the statement for Hodge classes defined by Weil. As in
SLN 900, an abelian variety with a space of Weil classes deforms smoothly in
characteristic zero to a power of an elliptic curve, on which all Hodge
classes are Lefschetz. The problem is that the family may not reduce
smoothly to $\mathbb{F}{}$. The final step will be to show that this gives enough classes for which the conjecture is true that the group fixing them is the correct one.

\section{ The case of Shimura varieties of Hodge type}

From now on I assume Conjecture A --- it seems to me essential to have such
a statement to obtain the canonical form of the conjecture Langlands and
Rapoport, even for Shimura varieties of PEL-type.

Then I can prove that, for a Shimura variety of Hodge type, there is a
canonical injection 
\begin{equation*}
LR(\mathbb{F}{})\rightarrow \Sh_{p}(\mathbb{F}{})
\end{equation*}

\noindent compatible with the actions of $G(\mathbb{A}_{f}^{p})$ and the
Frobenius automorphism.

The main difficulty in proving this statement involves the lattices in the $%
p $-cohomology. In characteristic zero, they lie in the $p$-adic \'{e}tale
cohomology, and in characteristic $p$, they lie in the crystalline
cohomology. Fortunately, the relation between these cohomologies is now
rather well understood, especially in the case of good reduction. The proof
of the statement uses theorems of Blasius and Wintenberger, Wintenberger,
and Fontaine and Messing.

I now need to assume another statement (Conjecture 0.1 of Milne 1995), which
is proved in a manuscript of Vasiu (.... Part 2A). As of writing, the proof
of Vasiu has not been checked. Appeal to Vasiu's paper can be avoided (I
think) if one extends Deligne's theorem on Tannakian categories (that any
two fibre functors are locally isomorphic) from Tannakian categories over
fields to Tannakian categories over Dedekind domains. (I have no idea
whether such an extension is possible, or even true, but it would be of
considerable interest if it is).

Assuming this, the map $LR(\mathbb{F}{})\rightarrow \Sh_{p}(\mathbb{F}{})$
is surjective if and only if the following conjecture holds:

\begin{conjecture}[B]
Zink's theorem holds for Shimura varieties of Hodge type.
\end{conjecture}

[[Restate in terms of Mumford-Tate groups. Equivalent statement that Hodge
classes on abelian varieties (not necessarily of CM-type) reduce to rational
Tate classes. Hence Conjecture B is implied by the Hodge conjecture for
abelian varieties.]]

\section{The case of Shimura varieties of abelian type.}

Happily, the extension from Hodge type to abelian type has been taken care
of by Pfau. Specifically, he shows that if a ``refined'' form of the
conjecture of Langlands and Rapoport holds for Shimura varieties of Hodge
type, then the same form of the conjecture holds for all Shimura varieties
of abelian type.

To state the refined form of the conjecture, he defines a map $LR(\mathbb{F}%
{})\rightarrow \pi _{0}(\Sh_{p})$. Here $\pi_{0}(*)$ denotes the set of connected components of $*$. The refined form of the conjecture then
states that there is a bijection $LR(\mathbb{F}{})\rightarrow \Sh_{p}(%
\mathbb{F}{})$ compatible with the maps to $\pi _{0}(\Sh_{p})$ (and the
actions of the Frobenius automorphism and $G(\mathbb{A}_{f}^{p})$).

Unfortunately, rather than a single well-defined map $LR(\mathbb{F}%
{})\rightarrow \pi _{0}(\Sh_{p})$ Pfau defines only a distinguished class of
maps.

I claim that one gets a canonical such map, almost for free. For a Shimura
variety $\Sh(G,X)$, let $T=G/G^{\text{der}}$ and let $\bar{X}=T(\mathbb{R}%
{})/Im(Z(\mathbb{R}{}))$ where $Z$ is the centre of $G$. The image of $Z(%
\mathbb{R}{})$ in $T(\mathbb{R}{})$ contains the identity component, and so $%
\bar{X}$ is finite. Define $\Sh(T,\bar{X})$ to be the system $\{T(\mathbb{Q}%
{})\backslash \bar{X}\times G(\mathbb{A}_{f})/K\}$ with $K$ running through
the compact open subgroups of $T(\mathbb{A}{}_{f})$. This is not a Shimura
variety in the sense of Deligne's original definition, but Pink has pointed
out that the study of the boundaries of Shimura varieties suggests that
Deligne's definition be extended to allow $X$ to be finite covering of a
conjugacy class of maps $\mathbb{C}{}^{\times }\rightarrow G(\mathbb{R}{})$.
For this extended definition, $\Sh(T,\bar{X})$ is a Shimura variety. Now
(under our continuing assumption that $G^{\text{der}}$ is simply connected), 
$\pi _{0}(\Sh(G,X))=\Sh(T,\bar{X})$.

The definition Langlands and Rapoport extends easily to give a set $LR(T,%
\bar{X})(\mathbb{F}{})$ and it is easy to prove the conjecture in this case:
there is a canonical bijection $LR(T,\bar{X})(\mathbb{F}{})\rightarrow
Sh_{p}(T,\bar{X})(\mathbb{F}{})$. The canonical Langlands-Rapoport
conjecture for $\Sh_{p}(G,X)$ will give a commutative diagram 
\begin{equation*}
\begin{array}{ccc}
LR(G,X)(\mathbb{F}{}) & \rightarrow  & \Sh_{p}(G,X)(\mathbb{F}{}) \\ 
\downarrow  &  & \downarrow  \\ 
LR(T,\bar{X})(\mathbb{F}{}) & \rightarrow  & \Sh_{p}(T,\bar{X})(\mathbb{F}%
{}).
\end{array}
\end{equation*}

\noindent Using that $\Sh_{p}(T,\bar{X})(\mathbb{F}{})=\pi _{0}(\Sh_{p}(G,X))
$, we see that this implies Pfau's refined form of the conjecture. Now,
Pfau's arguments show that the canonical Langlands-Rapoport conjecture for
Shimura varieties of Hodge type implies the same conjecture for Shimura
varieties of abelian type.

\section*{References}

Grothendieck, A., Standard conjectures on algebraic cycles. 1969 Algebraic
Geometry (Internat. Colloq., Tata Inst. Fund. Res., Bombay, 1968) pp.
193--199 Oxford Univ. Press, London.

\end{document}